\documentclass[11pt]{amsart}

\usepackage[all]{xypic}

\usepackage{latexsym}

\usepackage{amssymb}
\usepackage{amsfonts}
\usepackage{amscd}
\usepackage{amsmath,amsthm}

\newtheorem{lemma}{Lemma}[section]

\newtheorem{lem}[lemma]{Lemma}
\newtheorem{prop}[lemma]{Proposition}
\newtheorem{thm}[lemma]{Theorem}
\newtheorem{cor}[lemma]{Corollary}

{

}

\theoremstyle{definition}

\theoremstyle{remark}

\numberwithin{equation}{section}

\newenvironment{pf}{\noindent{\bf Proof.}}{\hfill $\square$\medskip}

\def\CC{{\mathbb C}}

\def\NN{{\mathbb N}}

\def\PP{{\mathbb P}}

\def\RR{{\mathbb R}}

\def\ZZ{{\mathbb Z}}

\def\0ol{{\bar 0}}
\def\1ol{{\bar 1}}
\def\2ol{{\bar 2}}
\def\ol2{{\bar 2}}
\def\3ol{{\bar 3}}
\def\4ol{{\bar 4}}
\def\5ol{{\bar 5}}
\def\6ol{{\bar 6}}
\def\7ol{{\bar 7}}
\def\8ol{{\bar 8}}
\def\9ol{{\bar 9}}

\def\bold0{{\bf 0}}
\def\bold1{{\bf 1}}
\def\bold2{{\bf 2}} 
\def\bold3{{\bf  3}}
\def\bold4{{\bf 4}}
\def\bold5{{\bf 5}}
\def\bold6{{\bf 6}}
\def\bold7{{\bf 7}}
\def\bold8{{\bf 8}}
\def\bold9{{\bf 9}}

\def\P2Skly{\PP^2_{Skly}}

\def\GL{\operatorname {GL}}

\def\pd{{\operatorname {\partial}}}

\def\Tr{\operatorname {Tr}}

\def\Aut{\operatorname{Aut}}

\def\det{\operatorname{det}}

\def\dim{\operatorname{dim}}

\def\Fract{\operatorname{Fract}}

\def\GKdim{\operatorname{GKdim}}

\def\id{\operatorname{id}}

\def\min{\operatorname{min}}

\def\pd{{\partial}}

\def\sup{\operatorname{sup}}

\def\ul1{\operatorname{\underline{1}}}

\def\l{\leftarrow}

\def\d{\downarrow}

\def\a{\alpha}
\def\b{\beta}

\def\d{\delta}

\def\l{\lambda}

\def\s{\sigma}

\def\D{\Delta}

\def\dirlim{\mathop{\vtop{\baselineskip -100pt\lineskip -1pt\lineskiplimit 0pt
\setbox0\hbox{lim}\copy0\hbox to \wd0{\rightarrowfill}}}\limits}
\def\invlim{\mathop{\vtop{\baselineskip -100pt\lineskip -1pt\lineskiplimit 0pt
\setbox0\hbox{lim}\copy0\hbox to \wd0{\leftarrowfill}}}\limits}

\def\I11{{1 \kern -0.8pt \! \mbox{l}}}
\def\mumu{{\mu\kern-4.2pt\mu}}
\def\bfmu{{\mu\kern-4.2pt\mu}}
\def\2slash{\backslash \! \backslash}

\def\boxtimes{\setbox0\hbox{$\Box$}\copy0\kern-\wd0\hbox{$\times$}}

\begin{document}

\title[Free subalgebras]{Free subalgebras of the skew polynomial rings $k[x,y][t;\s]$ and $k[x^{\pm 1},y^{\pm 1}][t;\s]$
}

\author{S. Paul Smith}
\address{ Department of Mathematics, Box 354350, Univ.  Washington, Seattle, WA 98195}
\email{smith@math.washington.edu}



\maketitle

\begin{abstract}
Let $k$ be a field, $R$ a commutative $k$-algebra, and $\s$ a $k$-algebra automorphism of $R$. The skew 
polynomial ring $R[t;\s]$ is generated by $R$ and an indeterminate $t$ subject to the relations $ta=\s(a)t$ for all $a \in R$. 
This paper shows that for certain $R$ and appropriate $\s$ there are elements $a,b \in R$ such that the subalgebra of 
$R[t;\s]$ generated by $at$ and $bt$ is a free algebra. For example, if $\s$ is an automorphism of the polynomial ring $k[x,y]$, then 
the subalgebra of $k[x,y][t;\s]$ generated by $xt$ and $yt$ is free if and only if $\s$ is not conjugate to an automorphism 
of the form $x \mapsto a x+p(y)$, $y \mapsto b y+c$, for any $a,b,c \in k$ and $p(y) \in k[y]$. Similarly, if $\s$ is an automorphism of $k[x^{\pm 1},y^{\pm 1}]$ of the form $\s(x)=x^ay^b$ and $\s(y)=x^cy^d$, then the subalgebra of $k[x^{\pm 1},y^{\pm 1}][t;\s]$ generated by $xt$ and $yt$ is free if the spectral radius of ${{a \; \; b} \choose {c \; \; d}}$ is $>2$; indeed, 
$k[x^{\pm 1},y^{\pm 1}][t;\s]$ contains a free subalgebra if and only if the spectral radius of ${{a \; \; b} \choose {c \; \; d}}$ is $>1$.
\end{abstract}

\section{Introduction}

\subsection{The main results} 
Let $k$ be a field. Let $R$ denote the commutative polynomial ring $k[x,y]$ or its localization 
$k[x^{\pm 1},y^{\pm 1}]$. Let $\s$ be a $k$-algebra automorphism of  $R$. The skew polynomial ring $R[t;\s]$ is 
the vector space $R \otimes_k k[t]$ with multiplication defined so that $R\otimes 1$ and $1 \otimes k[t]$ are subalgebras and 
$t^nr=\s^n(r)t^n$ for all $r \in R$ and $n \ge 0$.  

We prove the following results.

\begin{thm}
\label{thm.B}
Let $\s$ be an automorphism of $k[x,y]$. Then $k[x,y][t;\s]$ contains a free subalgebra if and only if $\s$ is not conjugate to an automorphism of the form $x \mapsto a x+p(y)$, $y \mapsto b y+c$, for any $a,b,c \in k$ and $p(y) \in k[y]$.
\end{thm}

 \begin{thm}
\label{thm.monom.automs.A}
 Let $M={{a \; \; b} \choose {c\;\;d}} \in \GL(2,\ZZ)$ and let $\s$ be the
 automorphism of $k[x^{\pm 1},y^{\pm 1}]$ defined by $\s(x)=x^ay^b$ and $\s(y)=x^cy^d$. 
Then $k[x^{\pm 1},y^{\pm 1}][t;\s]$ contains a free subalgebra if and only if the spectral radius of $M$ is $>1$. 
In that case, $k\{xt^{2n},yt^{2n}\}$ is a free subalgebra of $k[x^{\pm 1},y^{\pm 1}][t;\s]$ for all $n$ such that 
the spectral radius of $M^{2n}$ is $\ge 2$.  
\end{thm}

Some of our results apply to skew polynomial rings over other commutative rings.

\subsection{Skew Laurent extensions} 
It is natural to view a skew polynomial ring $R[t;\s]$ as a subalgebra of a slightly larger algebra.

Let $R$ be a commutative $k$-algebra and $\s \in \Aut_k(R)$. The {\sf skew Laurent extension} 
$$
R[t^{\pm 1};\s] := \bigoplus_{n \in \ZZ} Rt^n
$$
is defined by declaring that it is the free left $R$-module with basis $\{t^i \; | \; i \in \ZZ\}$
and multiplication defined to be the $k$-linear extension of $(ft^i)(gt^j):=f\s^i(g)t^{i+j}$ for all $f,g \in R$ and all $i,j \in \ZZ$. 
When $\s$ is the identity this is the ordinary ring of Laurent polynomials $R[t^{\pm 1}]$. 

We make $R[t^{\pm 1};\s]$ a graded ring by setting $\deg(t)=1$ and $\deg(R)=0$.

 \subsection{}
We consider graded subalgebras 
$$
k\{at,bt\}  := \bigoplus_{n=0}^\infty (Vt)^n \; \subset \; R[t;\s]
$$
where $V:=ka+kb$ is a 2-dimensional subspace of $R$. The degree-$(n+1)$ component of  $k\{at,bt\}$ is 
$$
V\s(V)\ldots\s^n(V)t^{n+1}
$$
so $k\{at,bt\}$ is a free algebra if and only if $\dim_k\big(V\s(V)\ldots\s^n(V)\big) = 2^{n+1}$ for all $n \ge 0$.

\subsection{Conventions/Notation}

Throughout this paper $k$ is a  field and $R$ a commutative $k$-algebra. We always assume that $R$ is an integral domain and write $K$ for its field of fractions. Whenever we say ``free (sub)algebra'' we will mean ``free (sub)algebra on $\ge 2$ variables''.

We write $k[x,y]$ for the polynomial ring on two variables and $k[x^{\pm 1},y^{\pm 1}]$ for its localization obtained by inverting $x$ and $y$. 

If $(G,+)$ is an abelian group we write $G^\times$ for $G-\{0\}$.

We write $\rho(M)$ for the spectral radius of a matrix $M$.

We write $|S|$ for the cardinality of a set $S$.

If $u$ and $v$ are elements in a $k$-algebra $A$ we will write $k\{u,v\}$ for the $k$-subalgebra they generate.

\subsection{Relation to other work} 

\subsubsection{Free subalgebras of division algebras} 
In 1983, Makar-Limanov discovered that the division ring of fractions of the ring of 
differential operators $\CC[x,\pd/\pd x]$ contains a free subalgebra \cite{ML}.  
Since then the question of which division algebras contain free subalgebras  
 has been of  considerable interest. 

A good account of recent progress on the question of which division algebras contain free subalgebras can be found in 
Bell and Rogalski's paper \cite{BR1}. They make further progress on this question in \cite{BR2}.

 \subsubsection{}
Despite the many results on free subalgebras of division algebras, 
there seem to be no known examples of free subalgebras of $R[t;\s]$ when $R$ is a finitely generated commutative 
$k$-algebra. Theorems \ref{thm.B} and \ref{thm.monom.automs.A} provide a host of such examples. Furthermore, since 
$R$ is a domain in those cases, $R[t;\s]$ has a division ring of fractions which then has a free subalgebra. 

Because the division algebra of fractions of an algebra $A$ is much larger than $A$ it is much easier to find free subalgebras
of division algebras. For example, if $\s$ is the automorphism of the polynomial ring $\CC[z]$ defined by $\s(z)=z+1$, then $\CC[z][t;\s]$ does not contain a free algebra\footnote{The ring $\CC[z][t;\s]$ is isomorphic to the enveloping algebra of the 2-dimensional non-abelian Lie algebra.} but its the division ring of fractions, $\Fract\big(\CC[z][t;\s]\big)$,  is isomorphic to 
$\Fract\big(\CC[x,\pd/\pd x]\big)$ so $\Fract\big(\CC[z][t;\s]\big)$ has a free subalgebra.

 \subsubsection{Gelfand-Kirillov dimension}
 
Gelfand-Kirillov dimension is a function $\GKdim:\{\hbox{$k$-algebras} \} \to \RR_{\ge 0}\sqcup\{\infty\}$. We
refer the reader to \cite{KL} for its definition. 
We will use two properties:
\begin{enumerate}
  \item[(a)]
  If $B$ is a free algebra on $\ge 2$ variables, then $\GKdim(B)=\infty$. 
  \item 
  If $B$ is a subalgebra of $A$, then $\GKdim(B) \le \GKdim(A)$.
\end{enumerate} 
Thus, an algebra of finite GK-dimension does not contain a free subalgebra.


\subsubsection{Exponential growth}

Let $A=A_0 \oplus A_1 \oplus \cdots$ be  an $\NN$-graded $k$-algebra. We say $A$  is {\sf connected} if $A_0=k$;
 {\sf locally finite} if $\dim_k(A_n)<\infty$ for all $n$; {\sf finitely graded} if it is 
locally finite and finitely generated as a $k$-algebra  \cite{RZ}.

Such an $A$ has {\sf exponential growth} if
$$
\lim\sup_{n \to \infty} \Big(\dim_k (A_n)\Big)^{\!1/n} >1.
$$

The ur-example of an algebra with exponential growth is the free algebra, $k\langle x_1,
\ldots,x_r \rangle$ on $r \ge 2$ variables each having degree 1: the dimension of its degree $n$ component is $r^n$.

\subsubsection{}
It is natural to ask if a locally finite graded algebra having exponential growth has a graded 
free subalgebra on $\ge 2$ variables. 

Golod and Shafarevich \cite{GS} showed the answer to that question is 
``no'' by constructing a finitely generated, connected, graded $k$-algebra $A$ of exponential growth such that every
element in $A_{\ge 1}$ is nilpotent. Their algebra is not finitely presented.
The question of whether every {\it finitely presented} locally finite graded algebra having exponential growth contains a graded  free subalgebra on $\ge 2$ variables remains open. 

\subsubsection{}
Given a specific algebra of exponential growth one can ask if it contains a 
free subalgebra on $\ge 2$ variables. This paper shows the answer is ``yes'' for various algebras of the 
form $R[t;\s]$.  The question can be sharpened: given $R$ and $\s$, and particular $a,b \in R$, and an integer $n$, when is the 
subalgebra $k\{at^n,bt^n\}$ free? Most of our results that prove the existence of a free 
subalgebra give explicit $a,b$ and $n$, such that $k\{at^n,bt^n\}$ is a free algebra.

\subsubsection{Relation to the results in \cite{SPS}}

The results in this paper complement those in \cite{SPS} where the following result is proved.

\begin{thm}
Let $K/k$ be a finitely generated extension field of transcendence degree two. Let $X$ be a smooth projective surface
such that $k(X)\cong K$.
\begin{enumerate}
  \item 
If $\s$ is an automorphism of $X$ and also denotes the induced automorphism of $k(X)$, 
then $k(X)[t^{\pm 1};\s]$ contains a free subalgebra if and only if the 
spectral radius of the automorphism of the N\'eron-Severi group induced by $\s$ is $>1$. 
  \item 
If $\s$ is an automorphism of $K$, then $K[t^{\pm 1};\s]$ contains a free subalgebra if and only if the 
dynamical degree of $\s$ is $>1$. 
\end{enumerate} 
\end{thm}

The methods in this paper are more accessible to algebraists than those in \cite{SPS}. 
Because the results in \cite{SPS} are more general they are not  as sharp as those in this paper. 
For example, suppose $\s$ is such that $k(x,y)[t;\s]$ contains a free subalgebra of the form $k\{at^n,bt^n\}$ for suitable
$a,b \in k(x,y)$; the results in this paper generally produce a smaller $n$ than those in \cite{SPS}. All the free algebras in 
\cite{SPS} are of the form $k\{at^n,bt^n\}$ where the divisor of zeroes of $a^{-1}b$ is very ample. In this paper, some of the 
free algebras are of the form $k\{at^n,\s(a)t^n\}$ and the question of whether  the divisor of zeroes of 
$a^{-1}\s(a)$ is very ample does not enter into the argument.

\subsection{Acknowledgements}
\label{sect.Ack}
Dan Rogalski's paper \cite{R} prompted my interest in free subalgebras. He told me that a test case was the 
subalgebra $k\{xt,yt\} \subset k[x^{\pm 1},y^{\pm 1}][t;\s]$ where $\s(x)=xy$ and $\s(y)=xy^2$. By Theorem \ref{thm.monom.automs.A}, that subalgebra is free. I thank George Bergman for showing me that the argument I used to 
answer Rogalski's question could be simplified, and improved to show that the subalgebra of $k[x^{\pm 1},y^{\pm 1}][t;\s]$ generated by $xt^n$ and $yt^n$ is free for all $n \ge 1$. The argument used to prove Theorem \ref{thm.mono.autom} 
is, in part, based on Bergman's ideas.

\section{Observations}
 
\subsection{}
\label{ss.conjug.class}
The question of whether $R[t;\s]$ contains a free subalgebra depends 
only on the conjugacy class of $\s$ because if $\s,\tau \in \Aut_k(R)$, there is a $k$-algebra isomorphism 
$\Phi:R[t;\s] \to R[s;\tau\s\tau^{-1}]$ given by  $\Phi(ft^j):=\tau(f)s^j$ for $f \in R$ and $j \in \NN$.

\subsection{}
\label{ss.replace.ab.by.1c}
Let $a$ and $b$ be non-zero elements of $K$. We will show that $k\{at,bt\} \cong k\{t,a^{-1}bt\}$.

 Define $a_0=1$. For each integer $m \ge 1$ define $a_m:=a\s(a)\cdots \s^{m-1}(a)$
and $a_{-m}:=\s^{-m}(a^{-1}) \cdots \s^{-1}(a^{-1}) 
=\s^{-m}(a_m)^{-1}$. 
The fact that $a_m\s^m(a_n)=a_{m+n}$ for all $m,n \in \ZZ$
ensures that the  $k$-linear extension of the map $\Psi:K[t^{\pm 1};\s] \to K[t^{\pm 1};\s]$ defined by 
$\Psi(gt^m):=a_mgt^m$ for $g \in K$ is a graded $k$-algebra automorphism. 
Since $\Psi(t)=at$ and $\Psi(a^{-1}bt)=bt$, $\Psi$ restricts to a graded $k$-algebra isomorphism $k\{t,a^{-1}bt\} \to k\{at,bt\}$. 

Thus, $K[t^{\pm 1};\s]$ contains a free subalgebra of the form $k\{at,bt\}$ for some $a,b \in K^\times$ if and only if it 
contains a free subalgebra of the form $k\{t,ct\}$ for some $c \in K^\times$. 

\subsection{}
\label{ss.gk=34}
If $\s=\id_R$, then $R[t;\s]=R[t]$, the polynomial ring of with coefficients in $R$. If $R$ is a finitely
generated $k$-algebra, then the Gelfand-Kirillov dimension of $R[t]$ is $1+\GKdim(R)$.

\subsection{}
\label{sect.center}
If $\s^n=\id_R$, then the subalgebra of $R[t;\s]$ generated by $t^{n}$ is the commutative polynomial ring over $R$ on the indeterminate $t^{n}$ and $R[t;\s]$ is a free $R[t^n]$-module with basis $\{1,t,\ldots,t^{n-1}\}$ whence 
$\GKdim \big( R[t;\s] \big) = 1+\GKdim(R)$. Furthermore, $R[t;\s]$ embeds in the ring of $n \times n$ matrices over $R[t^n]$.
A matrix ring over a commutative ring never contains a free subalgebra on $\ge 2$ variables.

\subsection{}
Let $\s$ be the automorphism of $k[x^{\pm 1},y^{\pm 1}]$   given
by $\s(x)=x$ and $\s(y)=xy$.  Then $k[x^{\pm 1},y^{\pm 1}][t^{\pm 1};\s]$ is isomorphic to the group algebra 
of the discrete Heisenberg group, the subgroup of $\GL(3,\ZZ)$ consisting of upper triangular matrices 
with 1's on the diagonal. As is well-known, the growth rate of this group is 4 so 
$\GKdim\big(k[x^{\pm 1},y^{\pm 1}][t^{\pm 1};\s]\big)=4$. 

\subsection{}
Let $E$ be an elliptic curve and $K=k(E)(z)$, the field of rational functions 
over the function field of $E$. Artin and Van den Bergh \cite[Ex. 5.19]{AV} 
showed there is an automorphism $\s$ of $K$ and a finitely graded subalgebra $B \subset K[t;\s]$ such that $\GKdim(B)=5$. 

\subsection{}
By (the proof of) \cite[Cor. 5.17]{AV}, if $\s$ is an automorphism of a smooth projective surface $X$ such that the induced automorphism of the Neron-Severi group has an eigenvalue that is not a root of unity, 
then $k(X)[t^{\pm 1};\s]$ has exponential growth. 

\subsection{}
In  \cite{R}, Rogalski proved that when $K/k$ is a finitely generated field extension of transcendence degree 
the Gelfand-Kirillov dimension  of a finitely generated $k$-subalgebra of $K[t^{\pm 1};\s]$ 
is either 3, 4, 5, or $\infty$. Rogalski's proof uses ideas and results from complex dynamics and algebraic geometry.

\subsection{}  
A noetherian locally finite $\NN$-graded $k$-algebra never has exponential growth \cite[Thm. 0.1]{SZ} so can't contain a free algebra.

\section{Valuations}

Let $\nu:K  \to \RR \sqcup\{\infty\}$ be a valuation such that $\nu(a)=0$ for all $a \in k^\times$.

We note that if $\nu(x) \ne \nu(y)$, then $\nu(x+y) = \min\{\nu(x),\nu(y)\}$.

If $S$ is a subset of $K$ we write $\nu(S):=\{ \nu(x) \; | \; x \in S\}$. 

\begin{lem}
\label{lem.valn}
Let $x_1,\ldots,x_n \in K^\times$. 
\begin{enumerate}
\item
If $\nu(x_i)>\nu(x_1)$ for all $i \ge 2$, then $x_1+x_2+\cdots+x_n \ne 0$.
  \item 
If $|\{\nu(x_1),\ldots,\nu(x_n)\}|=n$, then 
 $\{x_1,\ldots,x_n\}$ is linearly independent over $k$.
  \item 
If $U$ is a $k$-subspace of $K$, then $\dim_k(U) \ge |\nu(U^\times)|$.
\end{enumerate}
 \end{lem}
 \begin{pf}
 (1)
 This is a small variation on \cite[Lemma 1, p.8]{Serre} where the result is proved for a discrete valuation.
 Multiplying the $x_i$s by $x_1^{-1}$ we can assume that $x_1=1$. Thus, $\nu(x_i)>0$ for all $i \ge 2$. Hence
 each $x_i$, $i \ge 2$, belongs to the maximal ideal of the valuation ring associated to $\nu$. 
 The result now follows from the fact that $1$ is not in this maximal ideal.
 
 (2)
 Suppose the statement is false. Then there is a non-empty subset $I\subset\{1,\ldots,n\}$ and $\l_i \in k^\times$ such that 
 $\sum_{i \in I} \l_ix_i=0$. Since $\nu(\l_ix_i)=\nu(x_i)$, there is $j \in I$ such that $\nu(\l_ix_i)>\nu(\l_jx_j)$ for all $i \in I-\{j\}$.
 It follows from (1) that  $\sum_{i \in I} \l_ix_i \ne 0$. This is a contradiction so (2) must be true.
 
(3)
This follows at once from (2).
 \end{pf}
 
 Lemma \ref{lem.valn}(3) will be used to obtain lower bounds on the dimensions of $k$-subspaces of $K$.

\begin{thm}
\label{thm.main}
Let $K/k$ be a field extension and $\s \in \Aut_k(K)$.  Let $a,b \in K$. Let $L$ be the smallest 
 $\s$-stable subfield of $K$ that contains $a$ and $b$. Suppose $\nu:L \to \RR\sqcup\{\infty\}$ is a valuation on 
$L/k$ such that $\infty \ne  \nu(a) \ne \nu(b)\ne \infty$. If there is a number $\b \in \RR$ such that $\nu(\s(z))=\b \nu(z)$ for all 
$z \in L^\times$, then $k\{at^n,bt^n\}$ is a free subalgebra of $K[t;\s]$ for all $n$ such that $|\b^n| \ge 2$. 
\end{thm}
\begin{pf} 
We note that $\nu\big(\s^n(z)\big)=\b^n \nu(z)$ for all $z \in K^\times$.
The subalgebra of $K[t;\s]$ generated by $at^n$ and $bt^n$ is isomorphic to the subalgebra of $K[t;\s^n]$ generated by $at$ and $bt$ so, after replacing $\s$ by $\s^n$ and $\b$ by $\b^n$, we can, and will, assume 
$|\b| \ge 2$ and $\nu(\s(z))=\b \nu(z)$ for all $z \in K$. 

Define $\D:=\{\nu(a),\nu(b)\}$, $V:=ka+kb$, and 
$$
V_{n}:=V\s(V)\s^2(V)\cdots\s^{n-1}(V). 
$$
Because $\nu(a) \ne \nu(b)$,  $\dim_k(V)=2$.
The degree-$n$ component of $k\{at,bt\}$ is $V_{n}t^{n}$. 
To prove this theorem it suffices to show that $\dim_k(V_n) = 2^n$. 

We will do this by showing that $|\nu(V_n^\times)| = 2^n$ for all $n \ge 1$ and then invoking 
Lemma \ref{lem.valn}(2) to conclude that $\dim_k(V_n) \ge |\nu(V_n^\times)|=2^n$.  

An induction argument on $n$ shows that $\{\nu(\s^n(V^\times))\} =\b^n\Delta$.
Hence $\nu(V_{n}^\times)=\D + \b \D +\dots +\b^{n-1}\D$. 
 By hypothesis, $|\D|=2$. 

 Suppose $|\nu(V_n^\times)| = 2^n$ but $|\nu(V^\times_{n+1})| < 2^{n+1}$. 
Then 
$$
e'_0+e'_1\b+\cdots +e'_n\b^n +\nu(a)\b^{n+1}=e_0+e_1\b+\cdots  + e_n\b^n +\nu(b) \b^{n+1}
$$
for some $e'_0,\ldots,e'_{n},e_0,\ldots,e_{n} \in \D$. 
Hence
\begin{equation}
\label{equality}
\big(\nu(a)-\nu(b)\big)\b^{n+1}=(e'_0-e_0)+(e'_1-e_1)\b+\cdots+(e'_{n}-e_{n})\b^{n}.
\end{equation}
The absolute value of the left-hand side of (\ref{equality}) is $|\nu(a)-\nu(b)||\b|^{n+1}$ and the absolute value of the right-hand side   is  
$$
\le \, |\nu(a) -\nu(b)| \frac{|\b|^{n+1}-1}{|\b|-1}.
$$
The equality in (\ref{equality})  therefore implies 
$$
|\b|^{n+1} \le \frac{|\b|^{n+1}-1}{|\b|-1}
$$
which is false because $|\b|\ge 2$. We deduce that $|\nu(V^\times_{n+1})| = 2^{n+1}$
and therefore  $k\{at,bt\}$ is free. 
\end{pf}

\section{Monomial automorphisms of $k[x^{\pm 1},y^{\pm 1}]$} 

Let $R=k[x^{\pm 1},y^{\pm 1}]$. Let $M={{a \; \; b} \choose {c\;\;d}} \in \GL(2,\ZZ)$.
The automorphism $\s:R \to R$ defined by
\begin{equation}
\label{eq.mono.autom}
\s(x):=x^ay^b \qquad \hbox{and} \qquad \s(y):=x^cy^d
\end{equation}
is called a {\sf monomial automorphism} of $K/k$.  

\subsection{The finite-order case}

As remarked in \S\ref{sect.center}, if $\s$ has finite order, then $R[t^{\pm 1};\s]$ is a finite module over its center for every 
commutative ring $R$ so neither  $R[t^{\pm 1};\s]$, nor its division ring of fractions in the case when $R$ is a domain, contains a 
free subalgebra on $\ge 2$ variables.

The order of a monomial automorphism $\s$ is equal to the order of $M$. Since we will obtain results showing that $k\{xt,yt\}$ is a free subalgebra of $k[x^{\pm 1},y^{\pm 1}][t;\s]$ for suitable $\s$ we will briefly note some relations satisfied by $xt$ and $yt$ when 
$\s$ has finite order.

\begin{lem}
\label{lem.root.of.I}
If $n \ge 1$ and $M^{n}=I$, then $k\{xt,yt\}$ is not free because 
$(xt)^n(yt)^n=(yt)^n(xt)^n$.  
\end{lem}
\begin{pf}
The hypothesis implies that $\s^{n}$ is the identity map. Therefore
$$ 
x\s(x) \ldots \s^{n-1}(x)y\s(y) \ldots \s^{n-1}(y)= y\s(y) \ldots \s^{n-1}(y)  x\s(x) \ldots \s^{n-1}(x).
$$
Hence $(xt)^n(yt)^n =(yt)^n(xt)^n$.
\end{pf}

  \begin{lem}
\label{lem.non.free.1}
The subalgebra $k\{xt,yt\} \subset k[x^{\pm 1},y^{\pm 1}][t;\s]$ is not free if
\begin{enumerate}
  \item 
  $\Tr(M)=0$, or 
  \item{}
  $\Tr(M)=1$ and $\det(M)=-1$, or
  \item{}
  $\Tr(M)= \det(M)=-1$, or
  \item 
 $\Tr(M) \in \{ \pm 1, \pm 2\}$ and $\det(M)=1$.
\end{enumerate}
\end{lem}
\begin{pf}
(1)
If $\Tr(M)=0$, then $M^4=I$ so 
$(xt)^4(yt)^4=(yt)^4(xt)^4$.
 
(2)
A calculation shows that $(xt)^2(yt)=(yt)^2(xt)$.  

(3)
A calculation shows that $(xt)(yt)^2=(yt)(xt)^2$. 

(4)
One can easily verify the following claims:
\begin{enumerate}
  \item[(a)] 
if $\Tr(M)=2$, then $(xt)(yt)^2(xt)=(yt)(xt)^2(yt)$;
  \item[(b)]  
if $\Tr(M)=-2$, then $(xt)^2(yt)^2=(yt)^2(xt)^2$;
\item[(c)] 
if  $\Tr(M)=1$, then $(xt)(yt)(xt)=(yt)(xt)(yt)$;
\item[(d)] 
if  $\Tr(M)=-1$, then $(xt)^3=(yt)^3$.
\end{enumerate}
In case (c), we also note that $M^6=I$ so   $(xt)^6(yt)^6=(yt)^6(xt)^6$.
\end{pf}

\subsection{}
A matrix $M \in \GL(2,\ZZ)$ has finite order if and only if $\rho(M)=1$, i.e.,  if and only if condition (1) or (4) in 
Lemma \ref{lem.non.free.1} holds. 
 If condition (2) or (3) in Lemma \ref{lem.non.free.1}  is satisfied, then $\rho(M)=\frac{1}{2}(1+\sqrt{5})$
and, conversely, if  $\rho(M)=\frac{1}{2}(1+\sqrt{5})$, either (2) or (3) holds. If $\rho(M) \ne 1$, then 
$\rho(M)$ is either $\frac{1}{2}(1+\sqrt{5})$ or $>2$.

The trace of $M$ is the sum of its eigenvalues so  $\rho(M)>1$ if $|\Tr(M)|>2$.

\begin{thm}
\label{thm.mono.autom}
Let $M={{a \; \; b} \choose {c\;\;d}} \in \GL(2,\ZZ)$ and define $\s \in \Aut_k(k[x^{\pm 1},y^{\pm 1}])$ by $\s(x)=x^ay^b$ and $\s(y)=x^cy^d$. 
\begin{enumerate}
  \item 
If $\rho(M) =1$, then $k[x^{\pm 1},y^{\pm 1}][t;\s]$ is a finite module over its center so does not contain a free algebra
on $\ge 2$ variables.
 \item 
  If $\rho(M) = \frac{1}{2}(1+\sqrt{5})$, then $k\{xt^2,yt^2\}$ is a free algebra but  $k\{xt,yt\}$ is not. 
  \item 
  If $\rho(M) > \frac{1}{2}(1+\sqrt{5})$, then $k\{xt,yt\}$ is a free algebra. 
\end{enumerate}
Thus, $k[x^{\pm 1},y^{\pm 1}][t;\s]$ contains a  free subalgebra if and only if $\rho(M)>1$.
\end{thm}
\begin{pf}
(1)
Suppose $M$ has a non-real eigenvalue. Since the eigenvalues are the zeroes of the 
characteristic polynomial $x^2-\Tr(M)x+\det(M)$, $\Tr(M)^2-4\det(M) <0$. Hence $\det(M)=1$ and $\Tr(M) \in \{0,\pm 1\}$. By Lemma \ref{lem.non.free.1}, $k\{xt,yt\}$ is not a free algebra.

Suppose $M$ has a single real eigenvalue.  
Then $\Tr(M)^2-4\det(M)=0$ whence $\Tr(M)=\pm 2$ and $\det(M)=1$; Lemma \ref{lem.non.free.1}(4) tells us that
$k\{xt,yt\}$ is not free. 

Suppose $M$ has two distinct real eigenvalues. Those eigenvalues must be $+1$ and $-1$
so $\Tr(M)^2-4\det(M)=4$ which implies $\Tr(M)=0$ and $\det(M)=-1$. Lemma \ref{lem.non.free.1}(1) tells us that
$k\{xt,yt\}$ is not free. 

(2)
The fact that $k\{xt,yt\}$ is not a free algebra is proved by parts (2) and (3) of Lemma  \ref{lem.non.free.1}. The fact that 
$k\{xt^2,yt^2\}$ is free follows from part (3) of the present theorem because $k\{xt^2,yt^2\}$ is isomorphic to the 
subalgebra of $k[x,y][t;\s^2]$ generated by $xt$ and $yt$. 

(3)
Let $\b$ be an eigenvalue for $M$ such that $|\b|\ge 2$.

If $bc=0$, then either $a$ or $d$ is equal to 
$\b$; but $ad=\pm 1$ so that cannot be the case. Hence $bc \ne 0$. 

Both ${{b} \choose {\b-a}}$ and ${{d-\b} \choose {c}}$ are $\b$-eigenvectors for $M$.
Thus $M$ has a $\b$-eigenvector of the form  ${1\choose \a}$.  Because $M{1\choose \a}=\b{1\choose \a}$,
$$
a+b\a = \b 
\qquad \hbox{and} \qquad
c+d\a=\a\b.
$$
If $\a=1$, then $a+b=c+d=\b$ whence $\det(M)=(a-c)\b = \pm 1$, contradicting the fact that $|\b| \ge 2$. 
Therefore $\a \ne 1$.

Let  $\nu$ be the valuation on $k[x^{\pm 1},y^{\pm 1}]$ defined by
$$
\nu\Big(\sum a_{ij}x^iy^j\Big):=\min\{i+ j\a \; | \; a_{ij} \ne 0\}.
$$
Let $\D:=\{\nu(x),\nu(y)\}=\{1,\a\}$.  A simple calculation shows that
$
\nu\big( \s\big(x^iy^j\big)\big) =\b\nu(x^iy^j)$ 
whence  $\{\nu(\s^n(x)),\nu(\s^n(y))\} = \b^n\D$. 

Let $A$ be the algebra generated by $xt$ and $yt$. Then $A_{n+1}=V_{n+1}t^{n+1}$ where $V_{n+1}$ 
is the linear span of 
$$
\big\{ x_0\s(x_1)\ldots \s^n(x_n) \; \big\vert \; x_i\in\{x,y\} \big\}.
$$
To prove the theorem we must show that $\dim_k(V_{n+1}) =2^{n+1}$. Obviously, $\dim_k(V_{n+1}) \le 2^{n+1}$.

It is clear that
\begin{equation}
\label{eq.nu.values}
\big\{ \nu\big(x_0\s(x_1)\ldots \s^n(x_n)\big) \; \big\vert \; x_i\in\{x,y\}\big\} = 
\Bigg\{ \sum_{i=0}^n  \d_i \; \Bigg\vert \; \d_i\in \b^i \D \Bigg\} .
\end{equation}
Since $\a \ne 1$, $\b^i \D$ has exactly two elements. 
Since $\nu$ is a valuation, $\dim_k(V_{n+1})$ is at least the number of elements in the right-hand side of (\ref{eq.nu.values}).
To complete the proof we show that the set  on the right-hand side of (\ref{eq.nu.values}) 
contains $2^{n+1}$ elements. To do that it suffices to prove the following claim.

\underline{Claim:}
If $\d_i,\d_i' \in \b^i\D$ and $\d_0+\cdots+\d_n= \d'_0+\cdots+\d'_n$, then $\d_i=\d_i'$ for all $i$. 

\underline{Proof:} We argue by induction on $n$. The claim is true for $n=0$. 

Suppose the claim is true for $n$ but false for $n+1$.
Then there are elements $\d_i,\d_i' \in \b^i\D$, $0 \le i \le n+1$, such that
\begin{equation}
\label{delta.sum}
\d_0+\cdots+\d_n+\d_{n+1}= \d'_0+\cdots+\d'_n+\d'_{n+1}
\end{equation}
and some $\d_i \ne \d'_i$. If $\d_{n+1}=\d_{n+1}'$ the induction hypothesis implies that $\d_i=\d'_i$ for all 
$0 \le i \le n$. That is not the case so $\d_{n+1} \ne \d_{n+1}'$.
 
Since $|\d_j'-\d_j| \le |\b^{j}||1-\a|$, the absolute value of the right-hand side is  
$$
\Big\vert  \sum_{j=0}^n(\d_j'-\d_j) \Big\vert \le  |1-\a|  \sum_{j=0}^n |\b^{j}| < |1-\a| |\b^{n+1}| = |\d_{n+1}-\d_{n+1}'|
$$
 where the strict inequality follows from the hypothesis that $|\b|\ge 2$. Therefore
$$
\d_{n+1}-\d_{n+1}'  \, \ne \, \sum_{j=0}^n(\d_j'-\d_j).
$$
This contradicts (\ref{delta.sum}) so we conclude that the claim must be true for $n+1$. The validity of the claim completes the proof of the theorem. 
\end{pf}

\subsection{}
The automorphism of $k[x^{\pm 1},y^{\pm 1}]$ corresponding to the matrix  
$M={{1 \; \; 1} \choose {1 \; \; 2}}$ is given by $\s(x)=xy$, $\s(y)=xy^2$. The spectral radius of 
$M$ is $\frac{1}{2}(3+\sqrt{5})$ so $k\{xt,yt\}$ is a free subalgebra of $k[x^{\pm 1},y^{\pm 1}][t;\s]$. 
Thus, the answer to the question Rogalski asked---see \S\ref{sect.Ack}---is {\it ``yes''}.  

G. Bergman noticed that $\s$ is the square of the automorphism  
$\tau(x)=y$ and $\tau(y)=xy$ and showed that the subalgebra $k\{xt,yt\}\subset k[x^{\pm 1},y^{\pm 1}][t;\tau]$ is not free because 
$(xt)^2(yt)=(yt)^2(xt)$. Although $k\{xt,yt\}$ is not free it has exponential growth.  
The automorphism $\tau$ corresponds to $M'={{0 \; \; 1} \choose {1 \; \; 1}}$, $M=(M')^2$, and 
$\rho(M')=\frac{1}{2}(1+\sqrt{5})>1$. 
The subalgebra $k\{xt^2,yt^2\}\subset k[x^{\pm 1},y^{\pm 1}][t;\tau]$ {\it is} 
free because it is isomorphic to the subalgebra $k\{xt,yt\}$ in the previous paragraph.

\section{Automorphisms of $k[x,y]$}
\label{sect.k[x,y]}

\subsection{}
An {\sf elementary automorphism} of $k[x,y]$ is an automorphism $\tau$ of the form $x \mapsto a x+p(y)$, 
$y \mapsto b y +c$, for some $a,b,c \in k$, $p(y) \in k[y]$. 

A {\sf H\'enon automorphism} of $k[x,y]$ is an automorphism $\tau$ of the form 
$\tau(x, y) = (p(x)-ay, x)$ with $\deg(p(x)) \ge 2$ and $a \in k^\times$. A composition of H\'enon automorphisms is called a 
{\sf H\'enon map}.

\begin{thm}
[Friedland-Milnor]
\cite[Thm.2.6]{FM}
An automorphism of $\CC[x,y]$ is conjugate  to either  an elementary automorphism or  a H\'enon map.
\end{thm}

In this section we show that $k[x,y][t;\s]$ contains a free subalgebra if and only if $\s$ is conjugate to
a H\'enon map.

The dynamical system $(\CC^2,\s:\CC^2 \to \CC^2)$ where $\s$ is a H\'enon map of the form $(x,y)\mapsto (1+y-ax^2,bx)$ has been intensively studied in the context of complex dynamics. 
There is a belief that the ``chaotic behavior'' of any complex dynamical 
system is already exhibited by a H\'enon map for suitable $a$ and $b$. 
Guedj and Sibony say ``it is clear that only the [H\'enon maps] are dynamically interesting''  \cite{GS2}.

\begin{thm}
\label{thm.SPS.MKS}
Let $\nu$ be a valuation on a field extension $K/k$. Let $\s \in \Aut_k(K)$.  
Suppose there is a real number $\b>1$ and an element $g \in K$ such that either
\begin{enumerate}
  \item 
  $\nu(\s^{m+1}(g)) \ge \b \nu(\s^{m}(g))>0$ for all $m \ge 0$ or
  \item 
 $\nu(\s^{m+1}(g)) \le \b \nu(\s^{m}(g))<0$ for all $m \ge 0$.
\end{enumerate}  
If $\b^n \ge 2$, then $gt^n$ and $\s^n(g)t^n$ generate a free subalgebra of $K[t;\s]$.
\end{thm}
\begin{pf}
Define $\tau=\s^n$. If case (1) holds, then  $\nu(\tau^{m+1}(g)) \ge 2 \nu(\tau^{m}(g))>0$ for all $m \ge 0$. 
If case (2) holds, then  $\nu(\tau^{m+1}(g)) \le 2 \nu(\tau^{m}(g))<0$ for all $m \ge 0$. 

The subalgebra of $K[t;\s]$ generated by $gt^n$ and $\s^n(g)t^n$ is isomorphic to the subalgebra of 
$K[t;\tau]$ generated by $gt$ and $\tau(g)t$ so, after replacing $\s$ by $\tau$ and $\b$ by $\b^n$, we can, and will, assume $\b =2$ and either (1) or (2) holds.  

The first part of the proof applies to both cases (1) and (2).

Let $A$ be the algebra generated by $gt$ and $\s(g)t$. Then $A_{n+1}=V_{n+1}t^{n+1}$ where $V_{n+1}$ 
is the linear span of 
$$
\big\{ x_0\s(x_1)\ldots \s^n(x_n) \; \big\vert \; x_i\in\{g,\s(g)\} \big\}.
$$
To prove the theorem we must show that $\dim_k(V_{n+1}) =2^{n+1}$. Obviously, $\dim_k(V_{n+1}) \le 2^{n+1}$.

Define $\D_{n+1}:=\big\{ \nu\big(x_0\s(x_1)\ldots \s^n(x_n)\big) \; \big\vert \; x_i\in\{g,\s(g)\} \big\}$. Then
$$
\D_{n+1}:=\Bigg\{ \sum_{i=0}^n  \nu\big(\s^i(x_i)\big) \; \Bigg\vert \; x_i\in\{g,\s(g)\} \Bigg\}.
$$

Since $\nu$ is a valuation, $\dim_k(V_{n+1}) \ge |\D_{n+1}|$. We will complete the proof by showing 
that $|\D_{n+1}|=2^{n+1}$. This is 
true for $n=0$ because $\D_1=\{\nu(g),\nu(\s(g))\}$. Suppose the result is true for $\D_n$. Since
$$
\D_{n+1}=\Big(\D_n+\nu(\s^n(g))\Big) \; \bigcup \; \Big(\D_n+\nu(\s^{n+1}(g))\Big) 
$$
it suffices to show that 
$$
\Big(\D_n+\nu(\s^n(g))\Big) \; \bigcap \; \Big(\D_n+\nu(\s^{n+1}(g))\Big) =\varnothing.
$$
Let $\d_i=\nu(\s^i(g))$. To prove the intersection is empty it suffices to show that $\d_{n+1}-\d_n \ne \d-\d'$
for all $\d,\d' \in \D_n$. 

Now we split the proof into two separate parts according to the two cases in the statement of the theorem. 

(1)
In this case $0<\d_0 \le \frac{1}{2} \d_1 \le \cdots \le \frac{1}{2^n}\d_n$. 
The largest element in $\D_n$ is $\d_1+\cdots+\d_n=\nu\big(\s(g)\cdots \s^n(g)\big)$
and the smallest is $\d_0+\cdots+\d_{n-1}=\nu\big(g\s(g)\cdots \s^{n-1}(g)\big)$. If $\d,\d' \in \D_n$, then
$$
\d-\d' \le (\d_1+\cdots+\d_n)-(\d_0+\cdots+\d_{n-1}) = \d_n-\d_0 
$$
which is strictly smaller that $\d_{n+1}-\d_n $. 

(2)
In this case $0>\d_0 \ge \frac{1}{2} \d_1 \ge \cdots \ge \frac{1}{2^n}\d_n$. 
The smallest element in $\D_n$ is $\d_1+\cdots+\d_n=\nu\big(\s(g)\cdots \s^n(g)\big)$
and the largest is $\d_0+\cdots+\d_{n-1}=\nu\big(g\s(g)\cdots \s^{n-1}(g)\big)$. If $\d,\d' \in \D_n$, then
$$
\d-\d' \ge (\d_1+\cdots+\d_n)-(\d_0+\cdots+\d_{n-1}) = \d_n-\d_0 
$$
which is strictly larger that $\d_{n+1}-\d_n $. 
\end{pf}

\begin{cor}
 \label{cor.SPS.MKS}
Suppose $B$ be a commutative $\NN$-graded $k$-algebra and an integral domain. Let $\s \in \Aut_k(B)$ and   
$g \in B$. If $$\deg(\s^{m+1}(g)) \ge 2 \deg(\s^{m}(g))>0$$ for all $m \ge 0$, then $k\{gt,\s(g)t\}$ is a free subalgebra
of $B[t;\s]$.
\end{cor}
\begin{pf}
Write $K$ for the field of fractions of $B$. Let $\nu$ be the unique valuation on $K/k$ such that
$\nu(ab^{-1})=\deg(b)-\deg(a)$ whenever $a,b \in B$. The hypothesis in the statement of the corollary
implies that $g$ satisfies condition (2) in Theorem \ref{thm.SPS.MKS}. 
\end{pf}

\begin{cor}
\label{cor.MKS}
Let $\s$ be a $k$-algebra automorphism of $k[x,y]$.
\begin{enumerate}
  \item 
  $k[x,y][t;\s]$ contains a free subalgebra  if and only if $\s$ is not conjugate to an elementary automorphism.
  \item
  If $\s$ is not conjugate to an elementary automorphism, then  
  \newline 
   $k\{\s^n(x)t,\s^{n+1}(x)t\}$ is a free subalgebra of 
  $k[x,y][t;\s]$ for all $n \ge 0$.  
  \item 
  If  $\s$ is conjugate to an elementary automorphism, then
  \newline 
   $\GKdim\!\big(k[x,y][t;\s]\big)=3$. 
\end{enumerate}
\end{cor}
\begin{pf}
(3)
This is surely well-known but we could not find an argument in the literature so give one here. 

Suppose $\s$ is conjugate to a elementary automorphism. As noted in  \S\ref{ss.conjug.class}, the isomorphism class of $k[x,y][t;\s]$ as a graded $k$-algebra depends only on the conjugacy class of $\s$ 
so we can, and will, assume that $\s(x) = a x+p(y)$ and $\s(y)= b y +c$, for some $a,b,c \in k$, $p(y) \in k[y]$. 

Suppose $\deg(p)=d$. Let $V=k+ky+\cdots+ky^d+kx$ and let $W=V+kt$. 
Since $1,x,y,t \in W$, we can measure the GK-dimension of $k[x,y][t;\s]$ by measuring the 
rate at which $\dim_k(W^n)$ grows. Since $\s(V) =V$, $tV=Vt$.  An induction argument shows that
$$
W^n=V^n+V^{n-1}t+V^{n-2}t^2+ \cdots + Vt^{n-1}+kt^n.
$$
Therefore $\dim_k(W^n)=\dim_k(D^n)$ where $D$ is the subspace of the commutative polynomial ring 
$k[X,Y,T]$ spanned by $\{1,X,T,Y,Y^2,\ldots,Y^d\}$. Thus, $\GKdim\big(k[x,y][t;\s]\big) =\GKdim\big(k[X,Y,T]\big)=3$.

(1)
If $\s$ is conjugate to an elementary automorphism $k[x,y][t;\s]$ does not contain a free algebra on $\ge 2$ variables
because its GK-dimension is 3.

Suppose $\s$ is not conjugate to an elementary automorphism.  By \cite[Cor. 9]{MKS}, 
the degree of $\s^{n+1}(x)$ is at least twice the degree of $\s^n(x)$. By Corollary \ref{cor.SPS.MKS}, 
$k\{\s^n(x)t,\s^{n+1}(x)t\}$ is a free subalgebra of $k[x,y][t;\s]$ for all $n \ge 0$.  This completes the proof of (1) and also proves (2).
\end{pf}

\subsection{}
Let $a,b \in \CC$ with $ab \ne 0$. Let $\s$ be the automorphism of $\CC[x,y]$ defined by 
\begin{equation}
\label{defn.Henon}
\s(x)=1+y-ax^2  \qquad \hbox{and} \qquad \s(y)=bx.
\end{equation}
By Corollary \ref{cor.MKS}, $\CC\{xt,\s(x)t\}$ is a free subalgebra of $\CC[x,y][t;\s]$.

\begin{prop}
If $\s$ is the automorphism of $\CC[x,y]$ given by (\ref{defn.Henon}), 
then $\CC\{xt,yt\}$ is a free subalgebra of $\CC[x,y][t;\s]$.
\end{prop}
\begin{pf}
For the duration of this proof we give $\CC[x,y]$ the grading determined by $\deg(x)=2$ and $\deg(y)=1$. 
Since 
$$
\deg\big(\s(x^iy^j)\big)=i\deg\big(\s(x)\big) + j\deg\big(\s(y)\big) = 4i+2j = 2\deg(x^iy^j),
$$
$\deg(\s(f)) \le 2\deg(f)$ for all $f \in \CC[x,y]$. Since $\deg(x)=2$, an induction argument shows that
$\deg\big(\s^n(x)\big) \le 2^{n+1}$. (As we will shortly show, $\deg\big(\s^n(x)\big) = 2^{n+1}$.)  

\underline{Claim:} the degree-$2^{n+1}$ component of $\s^n(x)$ is a non-zero scalar multiple of $x^{2^{n+1}}$.
\underline{Proof:} The claim is true for $n=0$ and $n=1$. Suppose the claim is true for $n$, i.e., there is $\l_n \in \CC^\times$ such that  $\s^n(x)= \l_n x^{2^n} + {\rm l.d.t}$ where l.d.t stands for {\it lower-degree terms}, a {\it term} being a non-zero scalar multiple of some $x^iy^j$. Hence 
$$
\s^{n+1}(x)= \l_n \s(x)^{2^n} + \s({\rm l.d.t}) = \l_n(1+y-ax^2)^{2^n} + \s({\rm l.d.t}).
$$
In the previous paragraph we observed that the degree of $\s({\rm l.d.t})$ is at most twice the degree
of ${\rm l.d.t}$ so the degree-$2^{n+1}$ component of $\s^n(x)$ is a non-zero scalar multiple of $x^{2^{n+1}}$. 
$\lozenge$

Thus, $\deg(\s^n(y))=\deg(\s^{n-1}(x)) = 2^{n-1} = 2\deg(\s^{n-1}(y))$ for all $n \ge 1$. By Corollary \ref{cor.SPS.MKS}, with $g=y$, $\CC\{xt,yt\}$ is a free algebra.
\end{pf}


\section{Big subalgebras}

\subsection{The definition}
{(Rogalski and Zhang \cite[p.435]{RZ}, \cite[Defn. 6.1]{R})}.
Let $R$ be a commutative $k$-algebra and $\s \in \Aut_k(R)$. 
A locally finite $\NN$-graded subalgebra $\oplus_{n=0}^{\infty} V_nt^n \subset R[t;\s]$, where each $V_n \subset R$, is a {\sf big subalgebra} of $R[t;\s]$ if some $V_n$ contains a unit of $R$, $u$ say,  
such that $\Fract(R)=\Fract(k[V_nu^{-1}])$.

 \begin{prop}
\cite[Cor. 2.4]{RZ} 
Let $R$ be a commutative $k$-algebra. If a single finitely graded subalgebra of 
$R[t^{\pm 1};\s]$ has exponential growth so does every finitely graded big subalgebra of $R[t^{\pm 1};\s]$. 
\end{prop}

Let $\s \in \Aut_k\big(k[x^{\pm 1},y^{\pm 1}]\big)$. The subalgebra of $k[x^{\pm 1},y^{\pm 1}][t;\s]$ generated by $\{t,xt,yt\}$ is  a big subalgebra for all $\s$. The next result shows that the subalgebra of $k[x^{\pm 1},y^{\pm 1}][t;\s]$ generated by 
$xt$ and $yt$ can be free without being a big subalgebra. 

\begin{lem}
Let $M={{a \; \; b} \choose {c \; \; d}} \in \GL(2,\ZZ)$ and let $\s$ be the automorphism of $k[x^{\pm 1},y^{\pm 1}]$ given
by $\s(x)=x^ay^b$ and $\s(y)=x^cy^d$. Assume $a+b \equiv c+d \, {\rm (mod \, 2)}$. If $\rho(M) >2$, then $k\{xt,yt\}$ is free but is not a big subalgebra of $k[x^{\pm 1},y^{\pm 1}][t;\s]$.
\end{lem}
\begin{pf}
Since $\rho(M) >2$, $k\{xt,yt\}$ is a free algebra by Theorem \ref{thm.mono.autom}.

Let $A=k\{xt,yt\}$ and write $A_n=V_nt^n$ where $V_n=V\s(V) \cdots \s^{n-1}(V)$ and $V=kx+ky$. 
 
\underline{Claim:}
If $v,v' \in \s^n(V^\times)$, then $\deg(v)\equiv \deg(v')\, {\rm (mod \, 2)}$.
\underline{Proof:} 
We will prove this by induction on $n$. All non-zero elements in $V$ have odd degree so the claim is true for 
$n=1$. Suppose the claim is true for $n$. Let $x^iy^j \in \s^n(V)$. Since 
$$
\deg\big(\s(x^iy^j)\big) =(a+b)i+(c+d)j \equiv (a+b)(i+j) \,  {\rm (mod \, 2)}
$$
$\deg\big(\s(x^iy^j)\big) \, {\rm (mod \, 2)}$ is the same for all $x^iy^j \in \s^n(V)$. Hence the claim is true for $n+1$.
$\lozenge$ 

It follows from the claim that $\deg(u) \equiv \deg(u') \, {\rm (mod \, 2)}$ for all $u,u' \in V_n^\times$. Therefore 
$V_nu^{-1} \subset k(x^2,xy,y^2)$ for all $u \in  V_n^\times$ and all $n\ge 0$. Thus, $A$ is not a big subalgebra of 
$k[x^{\pm 1},y^{\pm 1}][t;\s]$.
\end{pf}

 \end{document}